\RequirePackage{fix-cm}
\documentclass{svjour3}                     % onecolumn (standard format)
\smartqed  % flush right qed marks, e.g. at end of proof

\usepackage[toc]{appendix}
\usepackage{geometry}
\usepackage{url}
\usepackage{cite}

\usepackage{color}
\usepackage{apptools}
\usepackage{chngcntr}%
\usepackage{graphicx}
\usepackage{amsmath}
\usepackage{amsfonts}
\usepackage{amssymb }
\usepackage{float}
\usepackage{multicol}
\PassOptionsToPackage{numbers}{natbib}
\usepackage{array}
\usepackage{tabularx}

\usepackage{scalerel,stackengine}
\def\apeqA{\SavedStyle\sim}
\def\apeq{\setstackgap{L}{\dimexpr.5pt+1.5\LMpt}\ensurestackMath{%
  \ThisStyle{\mathrel{\Centerstack{{\apeqA} {\apeqA} {\apeqA}}}}}}

\usepackage{latexsym}
\usepackage[sort&compress]{natbib}
\AtAppendix{\counterwithin{lemma}{section}}

\DeclareMathOperator{\sgn}{sgn}
\DeclareMathOperator{\Sgn}{Sgn}

\def\R{\mathbb R}

\def\({\left(}
\def\){\right)}
\def\<{\langle}
\def\>{\rangle}

 \journalname{Nonlinear Dynamics}
\makeatletter
\everymath{\displaystyle}

\begin{document}
\title{Complex dynamics, hidden attractors and continuous approximation of a fractional-order hyperchaotic PWC system}

%\subtitle{Hidden chaotic attractors in fractional-order systems}

\titlerunning{Hidden chaotic attractors in fractional-order systems}
% if too long for running head

\author{Marius-F. Danca \and M. Fe\v{c}kan \and Nikolay V. Kuznetsov \and Guanrong Chen}

\institute{
Marius-F. Danca  \and Corresponding author\label{corrauthor}\at
Department of Mathematics and Computer Science\\
Avram Iancu University, 400380 Cluj-Napoca, Romania \\
and \\
Romanian Institute for Science and Technology\\
400487 Cluj-Napoca, Romania \\
\email{danca@rist.ro}
\and
M. Fe\v{c}kan \at
Department of Mathematical Analysis and Numerical Mathematics, \\
Faculty of Mathematics, Physics and Informatics,\\
Comenius University in Bratislava, Mlynsk\'a dolina, 842 48 Bratislava, Slovak Republic,\\
and\\
Mathematical Institute, Slovak Academy of Sciences, \\
\v Stef\' anikova 49, 814 73 Bratislava, Slovak Republic \\
\email{Michal.Feckan@fmph.uniba.sk}
\and
Nikolay V. Kuznetsov\at
Department of Applied Cybernetics, \\Saint-Petersburg State University, Russia\\
and\\
Department of Mathematical Information Technology, \\University of Jyv\"{a}skyl\"{a}, Finland\\
\email: {nkuznetsov239@gmail.com}
\and
Guanrong Chen \at
Department of Electronic Engineering, \\City University of Hong Kong, Hong Kong, China \\\email{eegchen@cityu.edu.hk}}

\date{Received: date / Accepted: date}

\maketitle

\begin{abstract}
In this paper, a continuous approximation to studying a class of PWC systems of fractional-order is presented. Some known results of set-valued analysis and differential inclusions are utilized. The example of a hyperchaotic PWC system of fractional order is analyzed. It is found that without equilibria, the system has hidden attractors.

\vspace{3mm}
\end{abstract}
\section{Introduction}

This paper studies a new class of Piece-Wise Continuous (PWC) Fractional-Order (FO) systems modeled by the following general Initial Value Problem (IVP):
\begin{equation}
\label{IVP0}
D_*^q x(t)=f(x(t)),~~
x(0)=x_0,~~~ t\in I=[0,\infty),
\end{equation}

\noindent where the PWC function $f:\mathbb{R}^n\rightarrow \mathbb{R}^n$ has the form of
\begin{equation}\label{f}
f(x(t))=g(x(t))+A(x(t))s(x(t)),
\end{equation}

\noindent in which  $q\in(0,1)$, $g:\mathbb{R}^n\rightarrow \mathbb{R}^n$ a scalar vector-valued function, at least continuous, with $s:\mathbb{R}^n\rightarrow \mathbb{R}^n$, $s(x)=(s_1(x_1),s_2(x_2),...,s_n(x_n))^T$ a vector-valued piece-wise function, with $s_i:\mathbb{R}\rightarrow \mathbb{R}$, $i=1,2,...,n$, real piece-wise constant functions (e.g. $\sgn$ functions), and $A_{n\times n}$ a square matrix of real functions. Let $\mathcal{M}$ be the discontinuity set. Moreover, let $D^q_*$ denote the Caputo differential operator of order $q$ with starting point $0$ \cite{capa}:

\[
D_*^q x(t)=\frac{1}{\Gamma(1-q)}\int_0^t (t-\tau)^{-q}x'(t)d\tau.
\]

One of the reasons to use Caputo differential operator is that it has a physically meaningful interpretation for the initial conditions just like in the integer-order problems (a unique condition $x(0)$ in the case of $q\in(0,1)$).

\begin{remark}
Fractional-order differential equations (FDEs) do not define dynamical systems in the usual sense: by denoting the solution of \eqref{IVP0} as $\Phi(t,x_0)$, one does not have the flow property $\Phi_s\circ\Phi_t=\Phi_{t+s}$ \cite{zho}. However, in this paper, by numerical calculation of solutions, the definition of an integer-order dynamical system is adopted, which states that if the underlying IVP admits a solution, the problem defines a dynamical system (\cite[Definition 2.1.2]{stu}).
\end{remark}

Because the systems modeled by the IVP (\ref{IVP0}) are autonomous, hereafter, unless otherwise indicated, the time variable will be dropped in writing.

For numerical integration of discontinuous ordinary differential equations (ODEs) of integer order, there exist dedicated numerical methods (see e.g. the survey \cite{surv}) and, there are two main strategies to approach discontinuous systems: one strategy for treating discontinuities is simply to ignore them (\emph{time stepping methods}) and to rely on a local error estimator such that the error remains acceptably small; the other strategy is to use a scalar \emph{event function} $h:\mathbb{R}^n\rightarrow \mathbb{R}$, which defines the discontinuity $\Sigma=\{x\in \mathbb{R}^n|h(x)=0\}$, to determine the intersection point as the new starting point for the numerical solution (\emph{event-driven methods}).

The following aspects of discontinuous dynamical systems should be mentioned:

A numerical method for discontinuous systems may become either inaccurate or inefficient, or both, in a region where discontinuities of the solution or its derivatives occur and the local truncation error analysis, which forms the basis of most step-size control techniques, fails if there is not sufficient (local) smoothness.

Several known numerical methods assume that trajectories will cross the discontinuity surface as they reach it, or this will never happen. But there will always be errors in finding discontinuities.

Actually, systems of PWC systems are mostly ideal, since switch-like functions like $\sgn$ are used, where the hysteresis or delay of the real switching operation is not considered, or the regularization represents a good approach in these cases.

Although there are numerical methods for FDEs (see e.g. \cite{kai,kai2,dr}) and also for DEs with discontinuous right-hand sides (see e.g. \cite{don,bro,lem} or the survey \cite{surv}), to the best of our knowledge, there are no numerical methods for FDEs with discontinuous right-hand sides.
Consequently, modeling continuously or smoothly the underlying systems is of real important.

To approximate the PWC problem \eqref{IVP0}-\eqref{f}, one has to regularize the right-hand side first using e.g. the Filippov regularization, transforming the discontinuous problem to a set-valued IVP, i.e. an FO ordinary Differential Inclusion (DI) of FO. Then, Cellina's Theorem ensures the existence of continuous approximations\footnote{Graphically, by \emph{approximation} one understands that not all graphic points of the PWC function to be approximated need to be located on the created figure, compared with \emph{interpolation} where all graphic points of the PWC function have to be located on the created figure.}.

In this way, the obtained continuous FO problem can be numerically integrated using e.g. the multi-step predictor-corrector Adams-Basforth-Moulton (ABM) scheme for FDEs.

Two kinds of continuous approximations are proposed and utilized in this paper: global approximation and local approximation.

As an example of the PWC system modeled by \eqref{IVP0}-\eqref{f}, consider a fractional variant of one PWC system, proposed in \cite{spr1} for the integer case\footnote{See also \cite{spp}, where several PWC and nonsmooth jerk systems are proposed.}, as follows:

\begin{equation}\label{sp}
\begin{array}{cl}
D_*^{q}x_1= & -x_{1}+x_{2}, \\
D_*^{q}x_2= & -x_{3}\sgn(x_1)+x_{4}, \\
D_*^{q}x_3= & |x_{1}|-a,\\
D_*^qx_4=&-bx_2,
\end{array}%
\end{equation}

\noindent where $a,b>0$ are real positive parameters. In this paper, let $a=1$, $b$ be the bifurcation parameter and, unless otherwise specified, $q=0.98$.

Comparing \eqref{sp}  with \eqref{IVP0}, one has

\begin{equation*}
g(x)=\left(
\begin{array}{c}
-x_{1}+x_2 \\
x_4\\
|x_1|\\
-bx_2%
\end{array}%
\right) ,~~~A=\left(
\begin{array}{cccc}
0 & 0&0&0 \\
-x_3&0&0&0 \\
0&0&0&0%
\end{array}%
\right) ,~\ ~s(x)=\left(
\begin{array}{c}
\sgn(x_{1}) \\
\sgn(x_{2})\\
\sgn(x_3)\\
\sgn(x_4)%
\end{array}%
\right),
\end{equation*}

\noindent and $\mathcal{M}=\{0,x_2,x_3,x_4\}$, with the graph of $z=x_3\sgn(x_1)$ for $x_4=0$ shown in Fig. \ref{fig0} (a). Thus, $\mathbb{R}^4$ is splitted by the discontinuity surface, $x_1=0$, to two open half spaces, $\Omega_\pm=\{x\in\R^4 : \pm x_1>0\}$.

Besides rich dynamics (similar to those presented in \cite{spr1} for the integer-order case), it will be shown that the FO setting reveals some new behavior. It will be shown numerically that the local approximation to the PWC system by a continuous one is preferable instead of the classical global sigmoid approximation. Finite-time local Lyapunov exponents are determined, and chaotic as well as hidden hyperchaotic attractors are found. Sliding motion is numerically investigated. Moreover, using existing results on the periodicity of solution of FDEs, it will be shown that such a system cannot admit stable cycles.

The paper is organized as follows. Section 2 presents the continuous approximation of PWC systems \eqref{IVP0}. In Section 3, the PWC system \eqref{sp} is approximated. In Section 4, the dynamics of the obtained continuous system are numerically investigated. In the Conclusion section the obtained results are enumerated, while the three appendices present some mathematical notions and known results utilized by the paper.

\section{Approximation of the PWC problem }

\subsection{Regularization}

The PWC IVP \eqref{IVP0} will be transformed to the following FO DI:

\begin{equation}\label{DIgen}
D_*^qx\in F(x),~~x(0)=x_0,~~\text{for} ~a.a. ~t\in I.
\end{equation}

The set-valued function $F:\mathbb{R}^n\rightrightarrows 2^{\mathbb{R}^n}$ can be defined in several ways.
A simple (convex) expression of $F$, obtained by the Filippov regularization, is given by
\begin{equation}\label{fila}
F(x)=\bigcap_{\varepsilon >0}\bigcap_{\mu(\mathcal{M})=0} \overline{conv}(f({z\in \mathbb{R}^n: |z-x|\leq\varepsilon}\backslash \mathcal{M})),
\end{equation}
where $F(x)$ is the convex hull of $f(x)$, $\mu$ the Lebesgue measure and $\varepsilon$ the radius of the ball centered at $x$. At those points where $f$ is continuous, $F(x)$ consists of one single point, which coincides with the value of $f$ at this point (i.e. $F(x)=\{f(x)\}$). At the points belonging to $\mathcal{M}$, $F(x)$ is given by  \eqref{fila} (\cite[p.85]{filov}; see also \cite{waz}).

Set-valued function $F$, given by \eqref{fila}, is upper semicontinuous (see Definition \ref{usc} in Appendix), with compact and convex values (see also Remark \ref{rem} in Appendix). Then, applying \cite[Theorem 3.2]{zho} (see also \cite{sayed}), one obtains the following result about the existence of \emph{generalized} (\emph{Filippov}) solutions (see Definition \ref{gensol} in Appendix) to the DI (\ref{DIgen}).

\begin{theorem}\label{tui}
For every $x_0$, the differential inclusion \eqref{DIgen} with $F$ defined by \eqref{fila} has a generalized solution on $[0,\infty)$.
\end{theorem}

\noindent Even after regularization, the DI \eqref{DIgen} admits generalized solutions, because single-valued IVPs offer numerical opportunities. The interest in this paper is to transform the set-valued problem \eqref{DIgen} to a single-valued continuous one.

In order to justify the use of the Filippov regularization to some physical systems, one must choose small values for $\varepsilon$, so that the motion of the physical system is arbitrarily close to a certain solution of the underlying DI (it tends to the solution, as $\varepsilon\rightarrow 0$). However, extremely small values of $\varepsilon$ can lead to large values of derivatives and, consequently, to stiff systems. Therefore, as one can see in the following sections, special attention has to be focussed on $\varepsilon$.

If the piece-wise-constant functions $s_i$ are $\sgn$, their set-valued forms obtained with Filippov regularization and denoted by ${\Sgn}:\mathbb{R}\rightrightarrows \mathbb{R}$ are defined as follows:
\[
{\Sgn}(x)=\left\{
\begin{array}{cc}
\{-1\}, & x<0, \\
\lbrack -1,1], & x=0, \\
\{+1\}, & x>0.%
\end{array}%
\right.
\]
Here, the conventional $\sgn(0)$ is replaced with the whole interval $[-1 , 1]$ ''connecting'' the points $-1$ and $+1$ (see Figs. \ref{caca1} (a),(b)).

By applying the Filippov regularization to the function $f$ defined by \eqref{f}, one arrives at the following problem:

\[
D_*^qx\in F(x),~~~ x(0)=x_0,~~~ \text{for a.a.}~~t\in I,
\]
\noindent where

\begin{equation}
\label{IVP1}
F(x):=g(x)+A(x)S(x),
\end{equation}

\noindent in which
\[
S(x)=(S_1(x_1),S_2(x_2),...,S_n(x_n))^T,
\]

\noindent which $S_i:\mathbb{R}\rightarrow \mathbb{R}$ being the set-valued variant of $s_i$, $i=1,2,...,n$ (${\Sgn}(x_i)$ in the case of $\sgn$ function).

\subsection{Continuous approximation}

The continuous approximation of the function $f$ defined by \eqref{f} is realized next, as described in \cite{dancus}. For brevity, only the most important steps are presented.

\noindent Let $\mathcal{C}_\varepsilon^0(\mathbb{R})$ be the class of real continuous approximations (see Definition \ref{app} in Appendix) $\widetilde{s}:\mathbb{R}\rightarrow\mathbb{R}$ of the set-valued function $F$, which satisfy

(i) $Graph(\widetilde{s})\subset Graph (B(F,\varepsilon))$;

(ii) For every $x\in \mathbb{R}$, $\widetilde{s}(x)$ belongs to the convex hull of the image of $F$.

\smallskip

\noindent Here,  $B(x,\varepsilon)$ is the disk of radius $\varepsilon$ centered at $x$. Fig. \ref{caca1} (c) presents the case of the set-valued function $F(x)=\Sgn(x)$.

The approximation of the set-valued function $S$, with the single-valued function $\tilde{s}$, will be denoted as
\[
\tilde{s}(x)\apeq S(x).
\]

The set-valued functions $S_i$, $i=1,2,...,n$, can be approximated due to the Approximate Theorem, called Cellina's Theorem (Theorem \ref{theo} in Appendix), which states that a set-valued function $F$, with closed graph and convex values (Remark \ref{rem} in Appendix), admits $\mathcal{C}_\varepsilon^0$ approximations.

By \emph{global approximation} (GA) of the PWC functions $S$ in \eqref{IVP1}, denoted $\tilde{s}$, one understands a function defined on the entire axis $\mathbb{R}$, while by \emph{local approximation} (LA), denoted $\tilde{S}_\varepsilon$, a function defined on some $\varepsilon$-neighborhood $[-\varepsilon,\varepsilon]$ of discontinuity, with $\varepsilon$ being a small positive number \cite{dancus}.

\begin{theorem} \label{tete}\cite{dancus}~ The PWC system \eqref{IVP0}-\eqref{f}, with $g$ continuous, can be locally or globally continuously approximated with the following problem:
\[
D_*^qx=\tilde{f}(x),~~~ x(0)=x_0,~~~t\in I,
\]
where $\tilde{f}$ is either a local or a global approximation of $f$.
If $g$ is smooth, then the approximation is smooth.
\end{theorem}

\vspace{3mm}
\noindent\emph{Global approximation}
\vspace{3mm}

Any single-valued function on $\mathbb{R}$, with the graph located in the $\varepsilon$-neighborhood, can be considered as a global approximation of $S$ by Cellina's Theorem (see the sketch in Fig. \ref{caca1} (c), or the Weiestrass Theorem \ref{wei}). However, some of the best candidates for $\widetilde{s}$ are the \emph{sigmoid} functions, which provide the required flexibility and to which the abruptness of the discontinuity can be easily smoothened. If $S(x)={\Sgn}(x)$, one of the mostly utilized sigmoid approximations is the following function $\widetilde{\sgn}$\footnote{Sigmoid functions include the ordinary arctangent such as $\frac{2}{\pi}\arctan\frac{x}{\varepsilon}$, the hyperbolic tangent used especially in modelling neural networks (see e.g. \cite{dancila1}), the logistic function, some algebraic functions like $\frac{x}{\sqrt{\epsilon+x^2}}$, and so on.}:
\begin{equation}\label{h_simplu}
{\Sgn}(x)\apeq\widetilde{\sgn}(x)=\frac{2}{1+e^{-\frac{x}{\delta}}}-1,
\end{equation}

\noindent where $\delta$ is a positive parameter which controls the slope of the curve in the neighborhood of the discontinuity $x=0$. In Fig. \ref{fig1}(a), $\widetilde{\sgn}$ is plotted for three distinct values of $\delta$.

The smallest $\varepsilon$ values, necessarily to embed $\widetilde{\sgn}$ within an $\varepsilon$-neighborhood of ${\Sgn}$ (as stated by Cellina's Theorem), depends proportionally on $\delta$. Note that for $x\neq 0$, $\widetilde{\sgn}$ is identical to the single-valued branches of ${\Sgn}$ (the horizontal lines $\pm1$) only asymptotically, as  $x\rightarrow \pm\infty$. For example, for $\delta=0.01$, at the point $x=0.06$, the difference is of order $10^{-3}$, even the two graphs look apparently identical at the underlying points $A$ or $B$ (Fig. \ref{fig1}(b)). To reduce the $\varepsilon$ value, e.g., to $10^{-4}$, $\delta$ should be of order $10^{-5}$. On the other hand, as one can see in Section \ref{apr-sprot}, lower values of $\delta$ does not necessary imply substantial reduced values of $\varepsilon$.

In order not to significantly affect the physical characteristics of the underlying system, it is desirable to approximate $S$ only on some tight $\varepsilon$-neighborhoods of the discontinuity $x=0$, not on the entire real axis, since the difference between $S$ and $\tilde s$ persists along the entire real axis $\mathbb{R}$.

\vspace{3mm}
\noindent\emph{Local approximation}
\vspace{3mm}

A better approximation is LA, $\tilde{s}_\varepsilon:[-\varepsilon,\varepsilon]\rightarrow \mathbb{R}$, where $\tilde{s}_\varepsilon$ is some continuous function satisfying the gluing conditions

\begin{equation}\label{cond}
\widetilde{s}_{\varepsilon}(\pm\varepsilon)=s(\pm\varepsilon).
\end{equation}
Obviously, $\tilde{s}_\varepsilon$ and $\tilde{s}$ are both $C_\varepsilon^0$ functions.

\noindent If $g$ is continuous, then for every $\varepsilon>0$, there exists an LA of $f$, $\tilde{f}_\varepsilon:\mathbb{R}^n\rightarrow \mathbb{R}^n$, such that \cite{dancus} (see Fig. \ref{fig222} (a))
\[
f(x)\approx \tilde{f}_\varepsilon(x)=g(x)+A(x) {\widetilde{s}_\varepsilon}(x),~~~ x\in[-\varepsilon,\varepsilon].
\]

\noindent Approximation $\tilde{s}_{\varepsilon }$ can also be continuously extended on $\mathbb{R}$, yielding a GA, $\tilde{s}$, as follows:
\begin{equation}\label{fumu}
\tilde{s}(x)=\left\{
\begin{array}{cc}
\tilde{s}_{\varepsilon}(x), & x\in [ -\varepsilon ,\varepsilon ], \\
s(x), & x\notin [-\varepsilon ,\varepsilon ].%
\end{array}%
\right.
\end{equation}

\noindent Among the simplest functions $\tilde{s}_{\varepsilon }$ which, compared to the exponential function in the GA \eqref{h_simplu}, are the cubic polynomials (called \emph{cubic splines}) $\widetilde{s}_{\varepsilon }:\mathbb{R}\rightarrow \mathbb{R}$ defined by
\[
\widetilde{s}_{\varepsilon}(x)=ax^3+bx^2+cx+d,~~ a,b,c,d \in \mathbb{R}, ~a\neq0.
\]

\begin{remark}\label{poli}
Polynomials have the advantage to be directly evaluated by computers with the four arithmetic operations of adding, subtracting, multiplication and division.
\end{remark}

\noindent By imposing, near the gluing conditions (\ref{cond}), the supplementary differentiability conditions at the boundary of the discontinuity neighborhood\footnote{Since $s_i$ are PW constant functions on $x\neq 0$, they are differentiable on $x\neq 0$.}
\[
\frac{d}{dx}\widetilde{s}_{\varepsilon}(\pm \varepsilon)=\frac{d}{dx}s(\pm \varepsilon),
\]

\noindent for the particular case of the ${\Sgn}$ function, the smooth LA function, denoted by $\widetilde{\sgn}_\varepsilon$, is
\begin{equation}\label{loco}
{\Sgn}(x)\apeq \widetilde{\sgn}_{\varepsilon}(x)=-\frac{1}{2\varepsilon^3}x^3+\frac{3}{2\varepsilon}x,~~~  x\in [-\varepsilon, \varepsilon ].
\end{equation}

\noindent Using (\ref{fumu}), ${\Sgn}$ can be continuosuly approximated on $\mathbb{R}$ by the following piece-wise function (see Fig. \ref{fig222} (b)):
\begin{equation}
{\Sgn}(x)\apeq \left\{\label{loc}
\begin{array}{cc}
\widetilde{\sgn}_{\varepsilon}(x), & x\in [-\varepsilon ,\varepsilon ],  \\
\pm1, \text{(or~~} \sgn(x)\text{)}, & x\notin [-\varepsilon ,\varepsilon].
\end{array}%
\right.
\end{equation}

\begin{remark}
Both GA \eqref{h_simplu} and LA \eqref{loco} are also smooth approximations.
\end{remark}

Concluding, in order to obtain numerical solutions to \eqref{DIgen}, the simplest way is to replace the discontinuous problem with the continuous (smooth) one, using one of the locally or globally approximations (see the sketch in Fig. \ref{schita}).

\section{Approximations and numerical integration}\label{apr-sprot}

Following the way the problem was transformed to DI, the PWC problem \eqref{sp} can be transformed to the following set-valued problem:
\begin{equation}\label{DIfinal}
D_*^qx\in F(x)=\left(
\begin{array}{c}
-x_{1}+x_{2} \\
-x_{3}{\Sgn}(x_{1})+x_4 \\
|x_{1}|-a \\
-bx_{2}%
\end{array}%
\right), ~~~x(0)=x_0,~~~ \text{for a.a.}~~ t\in I,
\end{equation}

\noindent or
\[
\begin{array}{l}
D_*^qx_1=-x_{1}+x_{2}, \\
D_*^qx_2\in -x_{3}{\Sgn}(x_{1})+x_4, \\
D_*^qx_3=|x_{1}|-a, \\
D_8^q=-bx_{2}.%
\end{array}%
\]

\begin{figure}
\begin{center}
\includegraphics[scale=0.5]{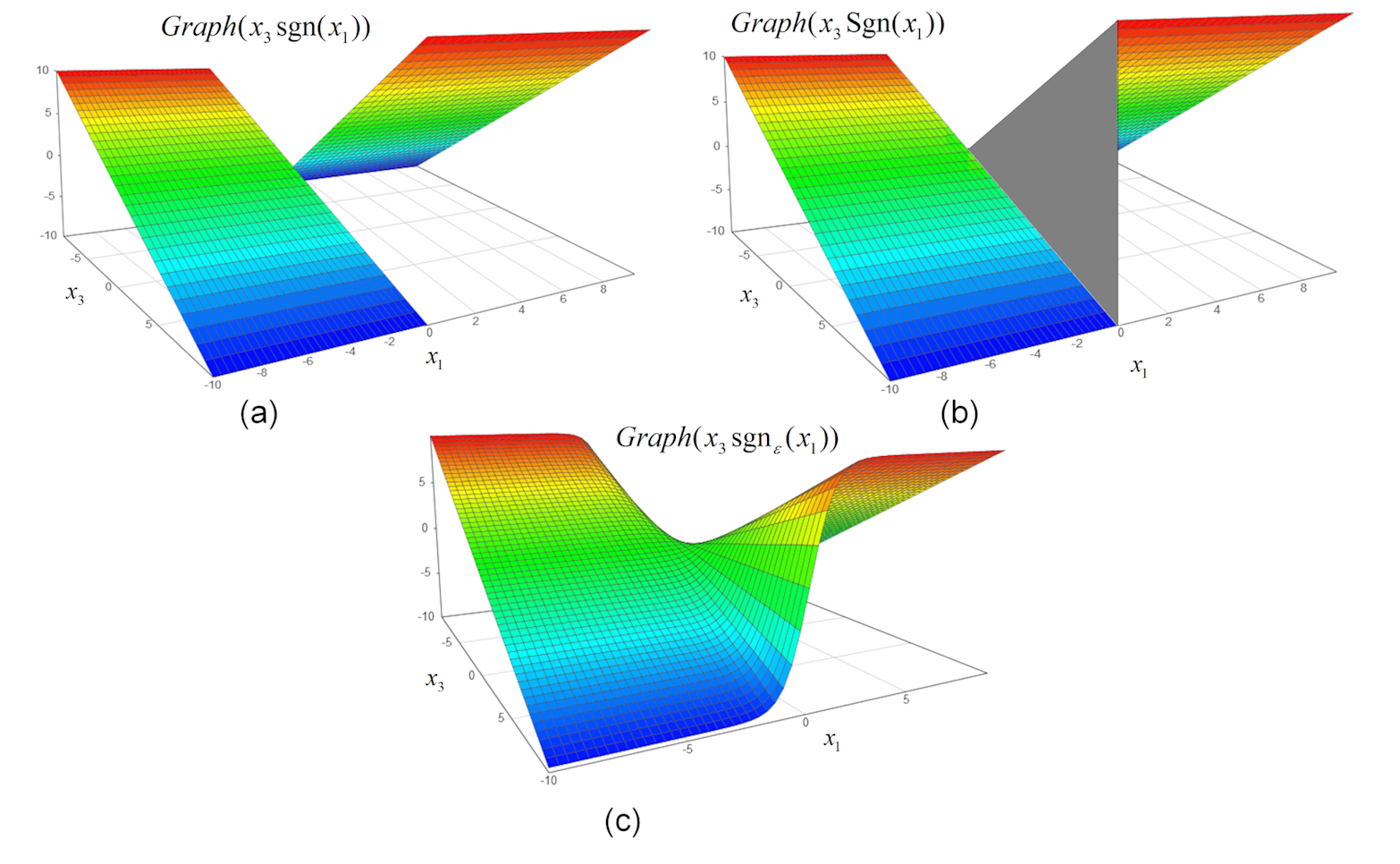}
\caption{Graph of the function $z=x_3\sgn(x_1)$. a) Before regularization. b) After regularization. c) After local continuous approximation for a large $\varepsilon$ value.}
\label{fig0}
\end{center}
\end{figure}

\begin{figure}
\begin{center}
\includegraphics[scale=0.5]{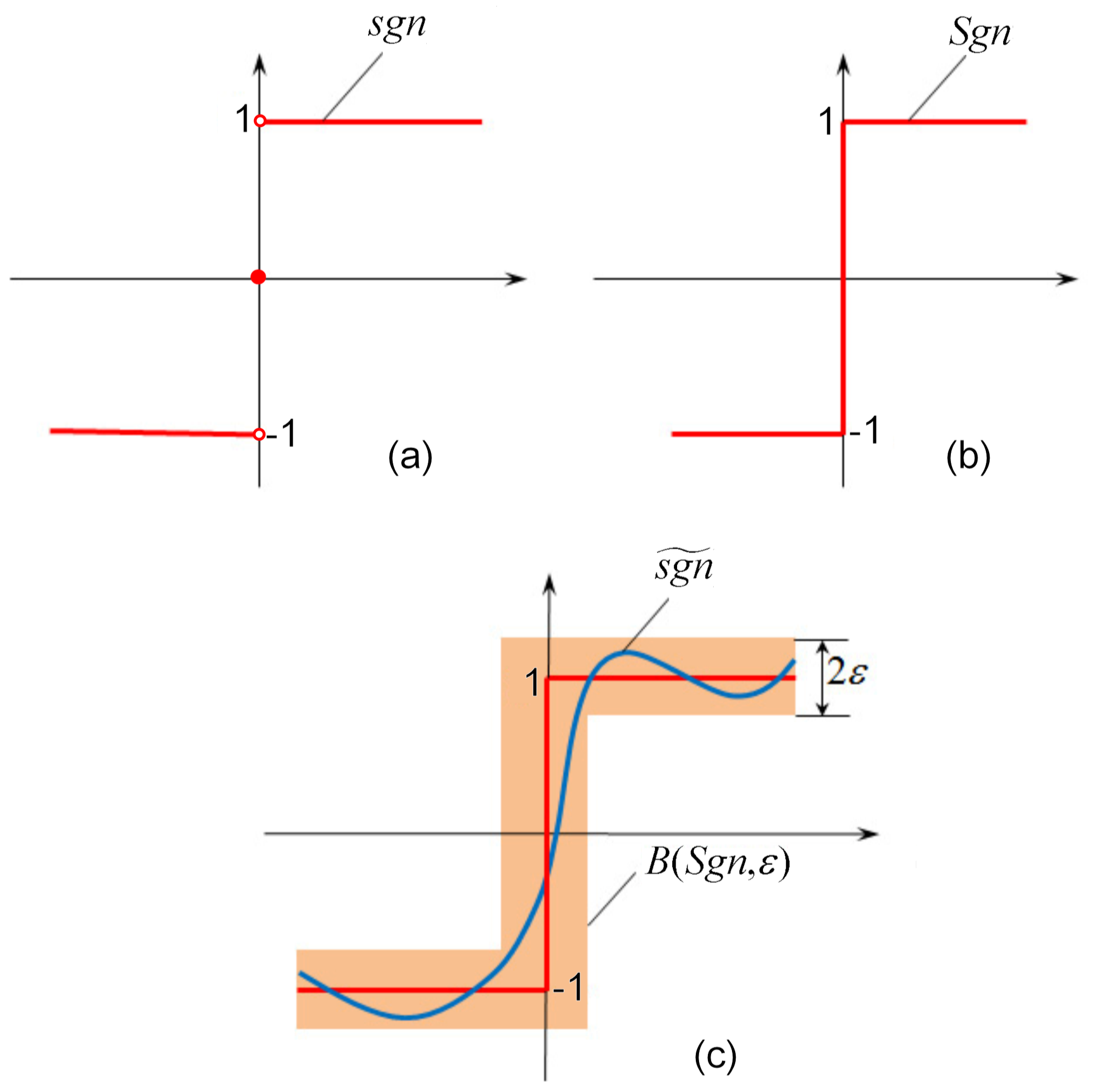}
\caption{(a) Graph of $\sgn$ function (red). (b) Graph of set-valued $\Sgn$ function (red). (c) Sketch of global approximation $\tilde{sgn}$ (blue) of $\Sgn$ function (red). $\Sgn$ neighborhood is plotted in light brown.}
\label{caca1}
\end{center}
\end{figure}

\begin{figure}
\begin{center}
\includegraphics[scale=0.5]{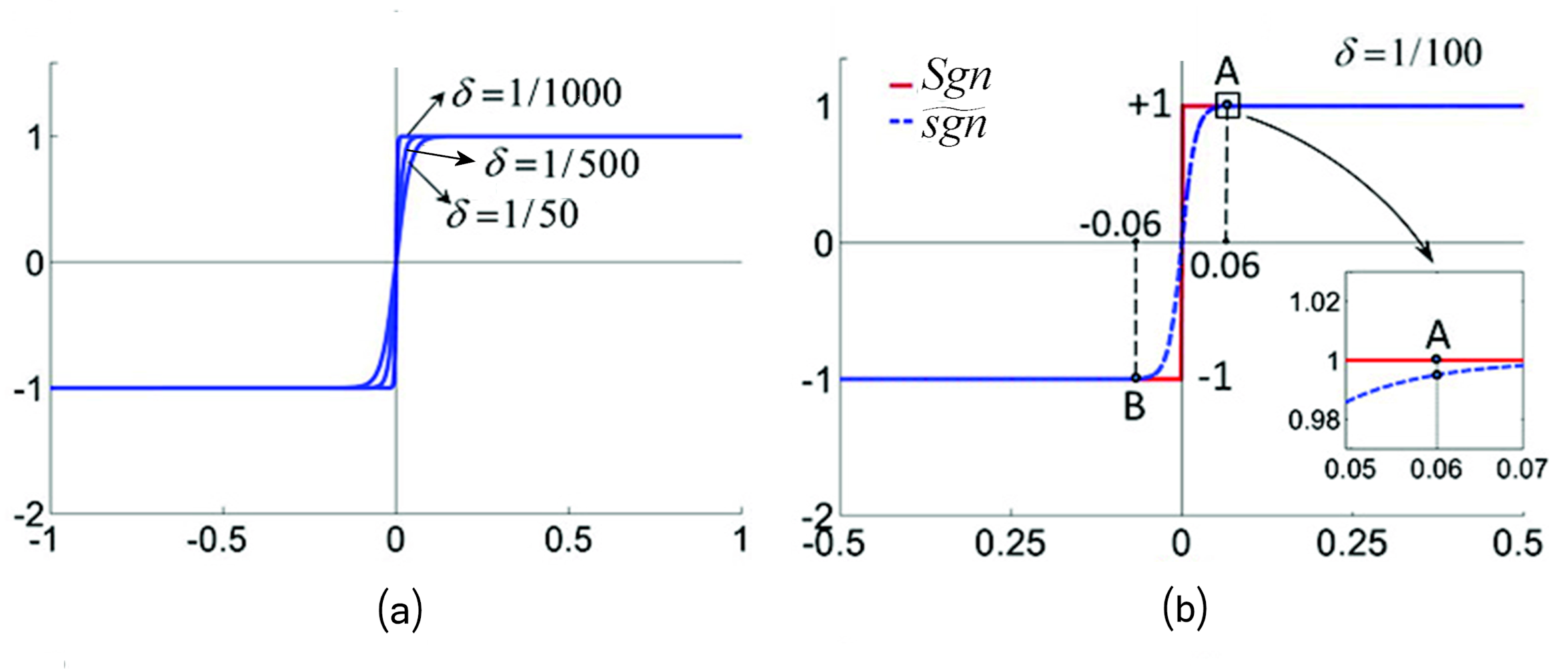}
\caption{a) Graph of sigmoid function $\widetilde{\sgn}$ \eqref{h_simplu} for three values of $\varepsilon$. b) Overplots of set-valued function $\Sgn$ and its approximate function $\widetilde{\sgn}$. Detail reveal the difference between the two curves.}
\label{fig1}
\end{center}
\end{figure}

\begin{figure}
\begin{center}
\includegraphics[scale=0.4]{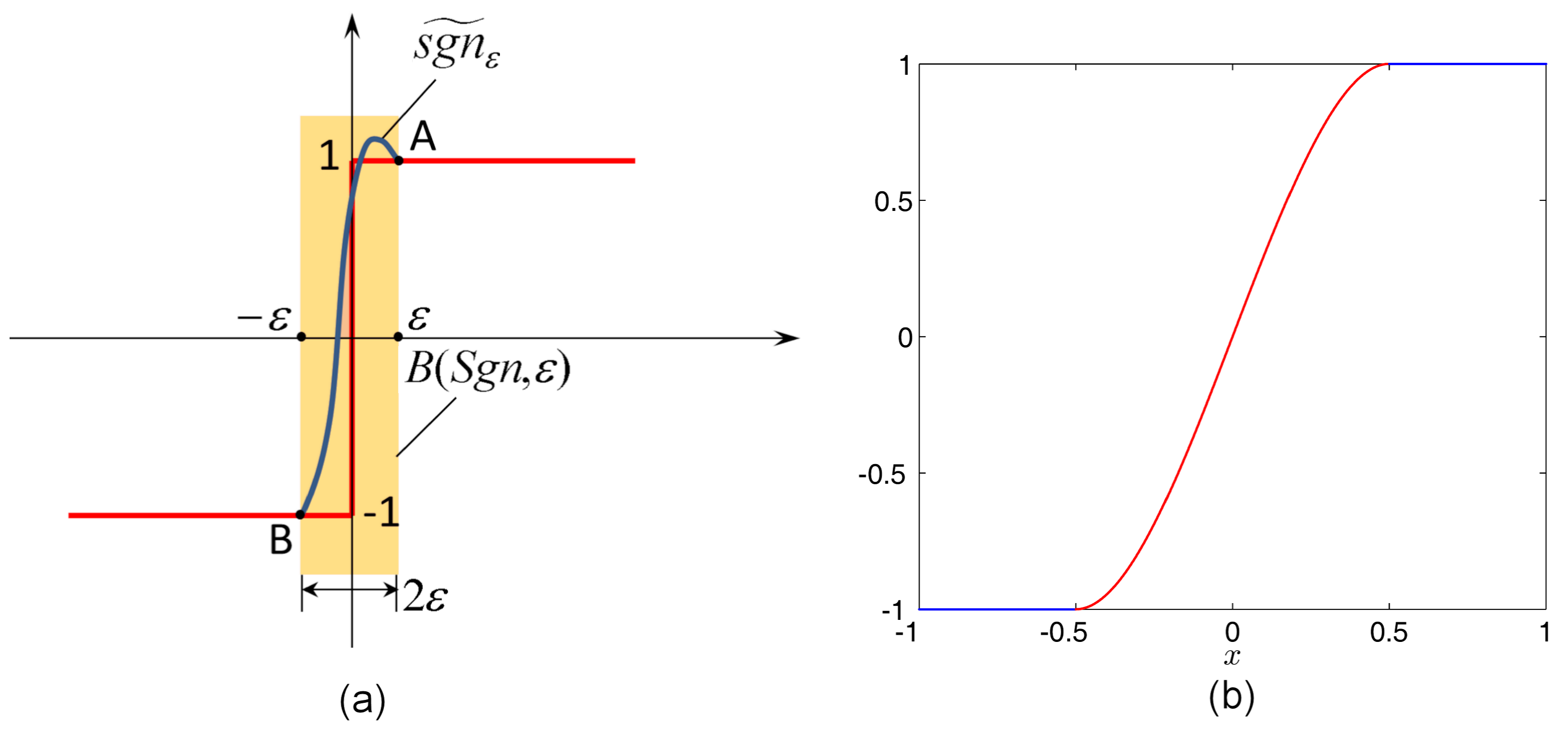}
\caption{(a) Sketch of local approximation of set-valued function $\Sgn$ (blue) within neighborhood $(-\varepsilon,\varepsilon)$. (b) Graph of cubic approximation $\widetilde{\sgn}_{\varepsilon}$ (red) in the neighborhood $(-\varepsilon,\varepsilon)$, with $\varepsilon=0.5$.}
\label{fig222}
\end{center}
\end{figure}

\begin{figure}
\begin{center}
\includegraphics[scale=0.5]{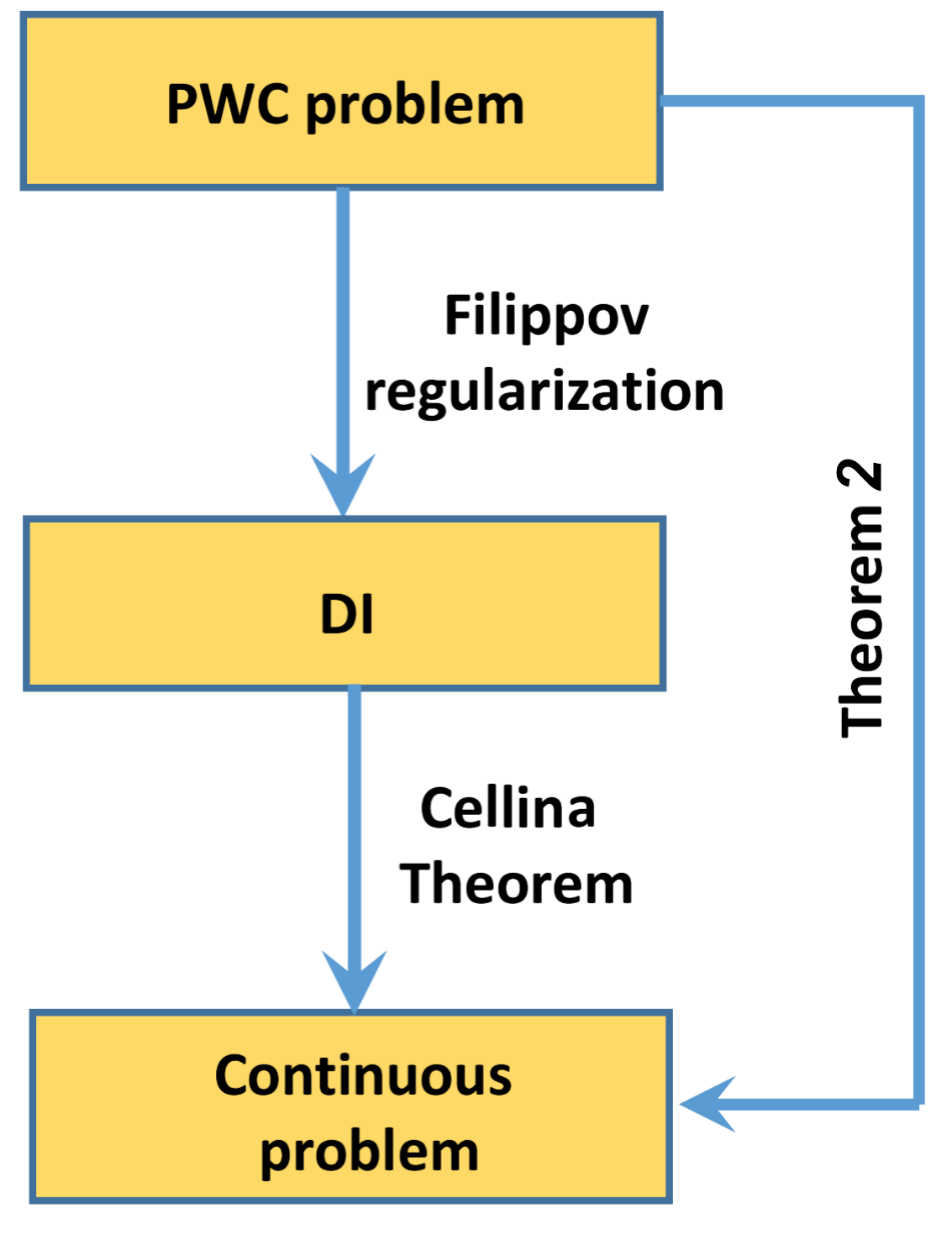}
\caption{Sketch of Theorem \ref{tete}.}
\label{schita}
\end{center}
\end{figure}

Note that only the second equation is a DI. For $x_4=0$, the underlying set-valued function is presented in Fig. \ref{fig0} (b).

By applying Theorem \ref{tete}, via relation \eqref{h_simplu}, GA leads to the following problem:
\begin{equation}\label{unu}
D_*^qx=\tilde{f}(x)=
\left(
\begin{array}{c}
-x_{1}+x_{2} \\
-x_{3}\left( \frac{2}{1+e^{\frac{-x_{1}}{\delta }}}-1\right)+x_4  \\
|x_{1}|-a \\
-bx_{2}.%
\end{array}%
\right),
\end{equation}

\noindent while with LA, \eqref{loc}, one has the following problem:

\begin{equation}\label{doi}
D_*^qx=\tilde{f}_\varepsilon(x)=\left(
\begin{array}{l}
-x_{1}+x_{2} \\
x_{4}-x_{3}\times\left\{
\begin{array}{l}
\widetilde{\sgn}_{\varepsilon }(x_{1}),~x_{1}\in \left[ -\varepsilon ,\varepsilon \right]
\\
\pm 1~\left( \text{or}~\sgn(x_{1})\right) ,~x_{1}\notin \left[ -\varepsilon
,\varepsilon \right]
\end{array}%
\right.  \\
|x_{1}|-a \\
-bx_{2}%
\end{array}%
\right).
\end{equation}

\begin{remark}\label{infini} Although both GA and LA are smooth, the approximating function $\tilde{f}$ is only continuous due to the modulus (third equation). However, the existence of solution to \eqref{DIfinal} is not affected by the non-smoothness (see e.g. \cite{di} for FO DIs).
\end{remark}

The graph of approximating surface of the second component on the right-hand side of the function $F$, with GA for a large value of $\delta$, $\delta=1$, is presented in Fig. \ref{fig0} (c). As can be seen, no significant graphical differences are indicated between both approximations.

The numerical integration of the approximate system is obtained in this paper with the predictor-corrector multi-step ABM method \cite{kai}, implemented using the Matlab subroutine FDE12.m \cite{roberto} with default tolerance $1E-6$, for fixed $q=0.98$ and, unless otherwise specified, integration step size $h=0.002$.

The following numerical analysis is performed for $b=1.25$ and $q=0.98$ with default double precision (roughly 15-16 decimal digits), sufficient for the current purpose.

\begin{figure}
\begin{center}
\includegraphics[scale=0.6]{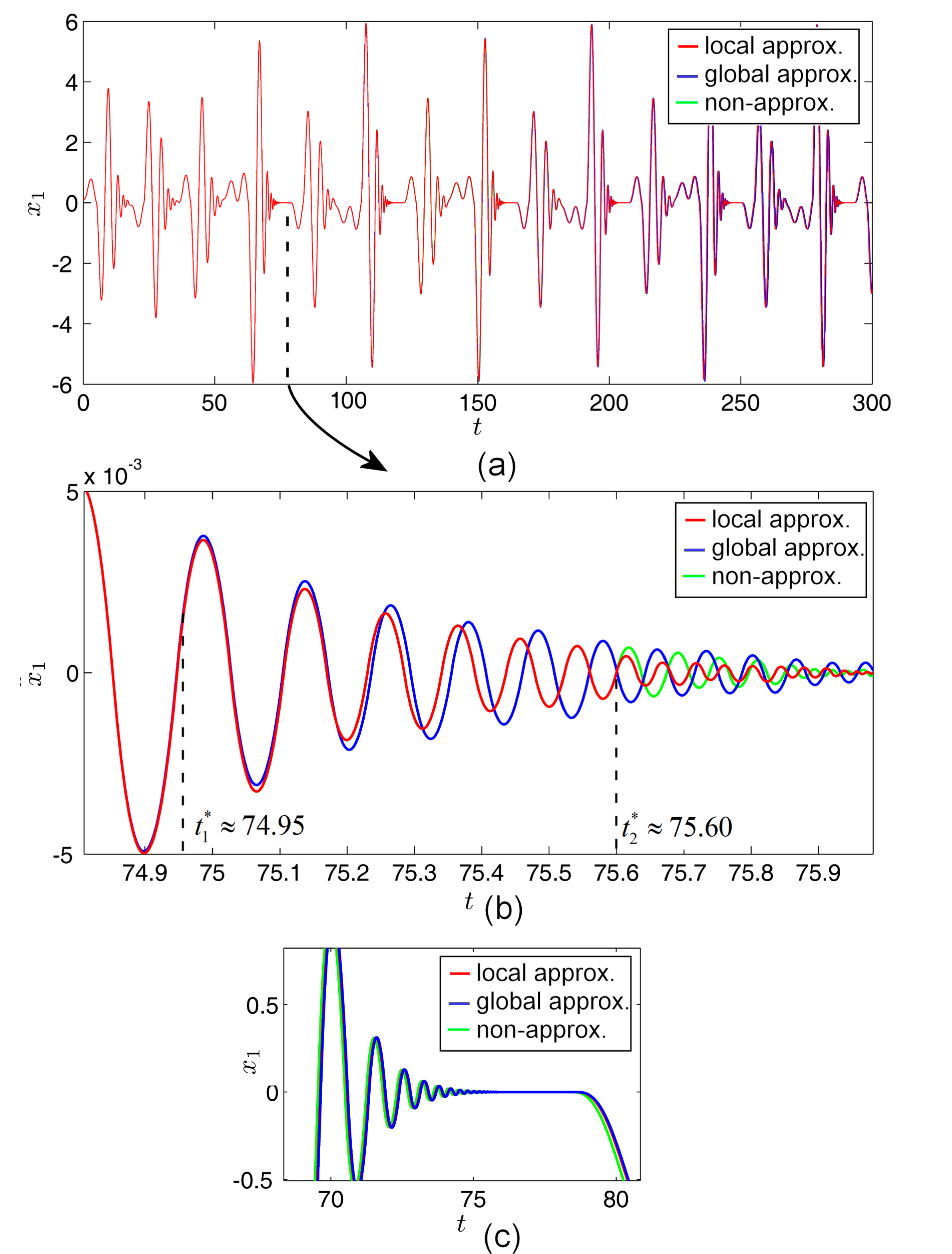}
\caption{(a) Overplotted time series $x_1$ of globally approximated system \eqref{unu} (blue plot), locally approximated system \eqref{doi} (red plot) and non-approximated system \eqref{sp} (green plot). (b) Detail. (c) Perfect identity of GA, LA and WA, along the horizontal line $x_1=0$, for step size $0.0002$.}
\label{pula}
\end{center}
\end{figure}

To compare the approximation results, Fig. \ref{pula} shows the time series from component $x_1$ for both GA and LA (with $\delta=1E-6$ and $\varepsilon=1E-6$, respectively), and also the time series obtained without approximation (WA) by using the ABM method on the non-approximated system, which are overplotted together.

The tests have been repeated with different initial conditions and time-step sizes.

Results are summarized as follows:
\begin{itemize}
 \item All the results (for GA, LA and WA) are accurate for all $t\in[0,t^*_1)$, with $t^*_1\approx 74.95$. Thus, for $t\in[0,t^*_1)$, the underlying numerical solutions coincide, because the differences between them are $0$;

  \item GA ``escapes'' early from this coincidence, near $t^*_1$ (Fig. \ref{pula} (b), while LA and WA remain further coinciding, till $t^*_2\approx75.60$;

  \item For decreasing values of $\delta$ and $\varepsilon$, $1E-10<\delta<1E-6$ and $1E-8<\varepsilon<1E-6$, the time superior limit slightly increases for both approximations ($t_{max}^*\nearrow76)$;

  \item For LA, with step size $0.002$ and $\varepsilon<1E-8$, the $\varepsilon$-neighborhood (necessary for the LA algorithm) can no longer be identified;

  \item For $\delta<1E-10$ in GA, numerically $\widetilde{\sgn}$ is the same as $\sgn$, therefore, there exists no more approximation, for which the code works further as FDE12 applied to the original WA system;

  \item With a smaller step-size of FDE12, no significant accuracy can increase. Therefore, up to $t\approx t_1^*$, both approximations are time-step independent (invariant) in the sense that the error between the two computed solutions, starting from the same initial conditions, remains zero, which are also close to WA;

  \item GA takes longer computational time, despite its simpler form, because it is calculated over the entire axis $x_1$, compared to LA which is calculated only on $[-\varepsilon,\varepsilon]$ (see also Remark \ref{poli}).
\end{itemize}

Note that, for a smaller step-size value (e.g. $h=0.0002$), in the crossing and sliding surface (i.e. $x_1=0$; see Subsection \ref{slid}), the two approximations are identical and also with the WA trajectory (Fig. \ref{pula} (c)).

Summarising, by comparing LA and GA and also with WA, it can be concluded that LA is better than GA in using with ABM for FDEs.

It is important to note that even FDE12.m can be applied directly to PWC problems, although the routine was not designed for this kind of problems.

\begin{remark}\label{timp}
\item []

The above results are in concordance with those obtained in \cite{lung1,lung2,lung3}.
Beyond numerical artifacts that might occur when numerically integrating a system of DEs, notions such as ``shadowing time'' and ``maximally effective computational time'' reveal that it is possible to have reliable numerical simulations only for some chaotic systems on a finite time interval (see e.g. \cite{lung1,lung2,lung3}). For example, in the classical Lorenz system, obtaining a precise solution for e.g. $t\in [0,100]$ represents a real challenge. Therefore, the case of FO systems is even more delicate.

On the other hand, larger intervals must be considered, so that possible phenomena like transient behaviours could be studied.

In this paper, after intensive numerical tests it is concluded that the obtained numerical results could be acceptable even on relatively larger intervals, as large as $[0,800]$.

\end{remark}

\vspace{3mm}

\section{Dynamics of the investigated FO PWC system}

Next, because of the above-mentioned advantages, the LA will be utilized.

\subsection{Sliding motion}\label{slid}

In order to check that there exist crossing, sliding and grazing phenomena in system \eqref{sp}, it is first considered without approximation.

 Denote the open half spaces by $\Omega_\pm=\{x\in\R^4 : \pm x_1>0\}$. For $x_0\in\Omega_\pm$, the switching time $t_s$ is given by the following equation (see Appendix \ref{apB}):
\[
\phi_{\pm}(t_s)=0,
\]
where
$$
\begin{gathered}
\phi_{\pm}(t)=e_1^T\left(E_q(t^qM_{\pm})x_0-at^qE_{q,q+1}(t^qM_{\pm})e_3)\right),
\end{gathered}
$$
with

$$
\begin{gathered}
e_1=\begin{pmatrix} 1\\ 0\\ 0\\ 0
\end{pmatrix},
e_3=\begin{pmatrix} 0\\ 0\\ 1\\ 0
\end{pmatrix} .
\end{gathered}
$$

Since $\phi_{\pm}(t)$ are analytic functions, they may only have isolated zeros. Consequently, crossing, sliding and grazing at $x(t_s)$ is possible, which are determined by the sign behaviour of functions $\phi_+(t)$ and $\phi_-(t)$ near $t_s$.

However, one cannot use the local behavior of the function $f$ given by \eqref{f}. Therefore, one does not know if the solution is really sliding on $x_1=0$, except the obtained numerical results. Moreover, one does not know how the solution is crossing the discontinuous surface. When a possible solution is crossing the discontinuous surface at some points, there is no theoretical proof for the existence of a solution.
Therefore, the graphical approach is adopted to study the trajectories of the continuous approximated system, to see what happens near the discontinuity surface $x_1=0$.

\begin{figure}
\begin{center}
\includegraphics[scale=0.5]{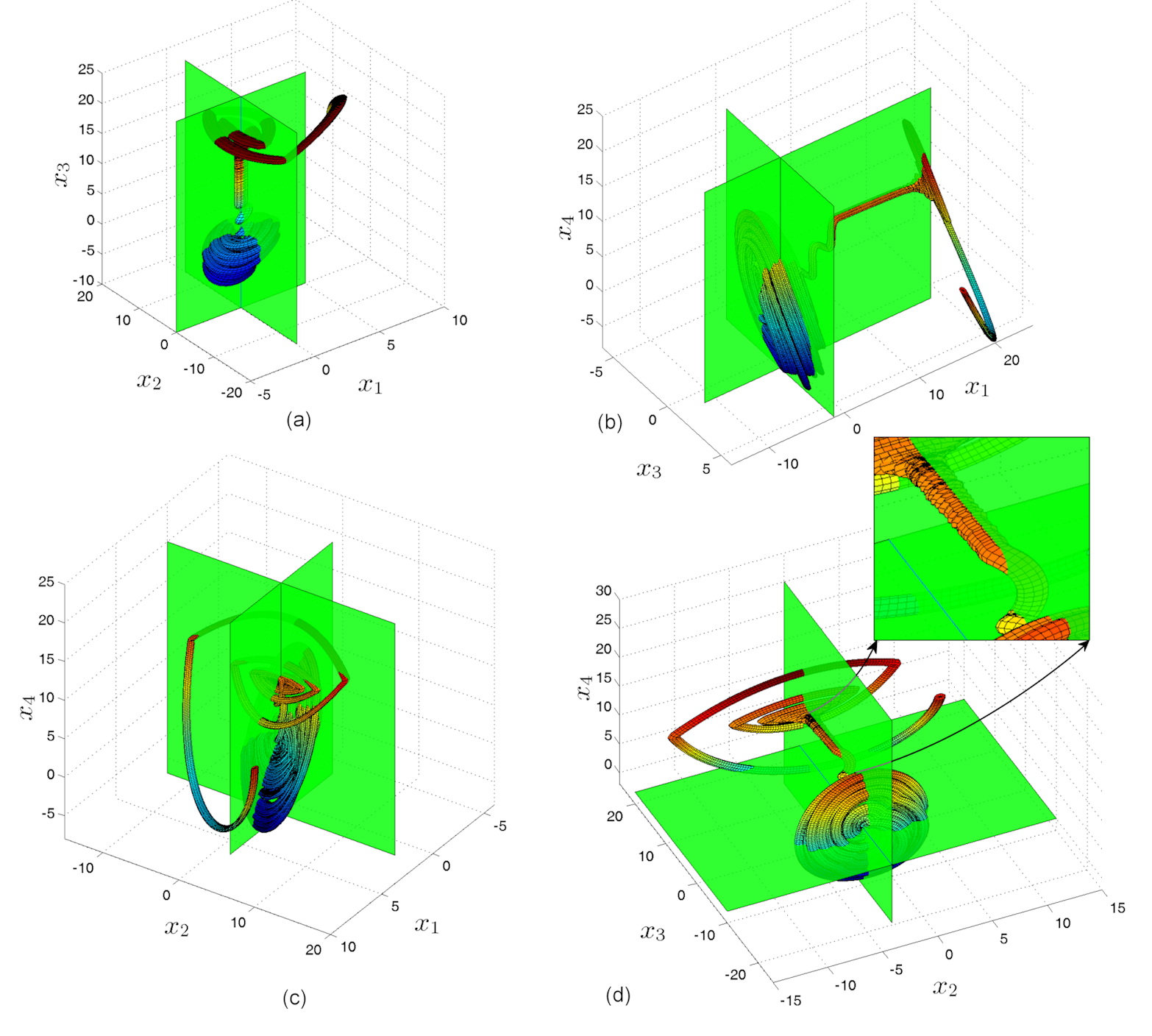}
\caption{3D phase projections of the hidden hyperchaotic attractor for $b=2.2$. (a) Space $(x_1,x_2,x_3)$. (b) Space $(x_1,x_3,x_4)$. (c) Space $(x_1,x_2,x_4)$. (d) Space $(x_2,x_3,x_4)$ and zoomed view. Possible sliding phenomenon can be observed in Fig. \ref{ciudat} (a), (b) and (d).}
\label{ciudat}
\end{center}
\end{figure}

As can be observed from Figs. \ref{ciudat} (a), (b) and (d), where the case of $b=2.2$ is considered, the trajectory seems to slide along the plane $x_1=0$. However, the tubular representations shown in Fig. \ref{ciudat} (especially the detail in Fig. \ref{ciudat} (d)) would indicate that the trajectory actually crosses the plane $x_1=0$ for several (but finite number of) times. It is remarked that this phenomenon happens for other values of $b$ too.

Note that, after approximation, the system still remains in class $C^0$ (see Remark \ref{infini}). Therefore, beside the ``smoothed'' corners caused by the approximated discontinuity, which appear along the plane $x_1=0$\footnote{As is well known, ODEs with $C^k$ class right-hand side have $C^{k+1}$ solutions.}, the trajectories have some other corners due to the modulus $|x_1|$ (see the 3D phase projections in Figs. \eqref{ciudat} (a), (c), (d)).

An interesting behavior has been found, especially when initial conditions are situated in a relatively larger distance from the origin of the phase space, where the trajectories scroll toward a direction parallel with the axis $x_3$ in the 3D phase projection $(x_1,x_2,x_3)$, or with some axis parallel to the axis $x_1$ in the plane $x_3=0$ in the 3D phase projection $(x_1,x_3,x_4)$ (see Fig. \ref{ciudat} (b)). This also happens, on the plane projection $(x_4,x_3)$, for the hidden hyperchaotic attractor corresponding to $b=0.5$ in Fig. \ref{fig222} (d).

These scrolls differ, for example, from the scrolls generated by nonlinear modulating functions in jerk systems (see e.g. \cite{jerk}, or \cite{spp}).

\subsection{Periodicity}\label{perioada}
\begin{figure}
\begin{center}
\includegraphics[scale=0.4]{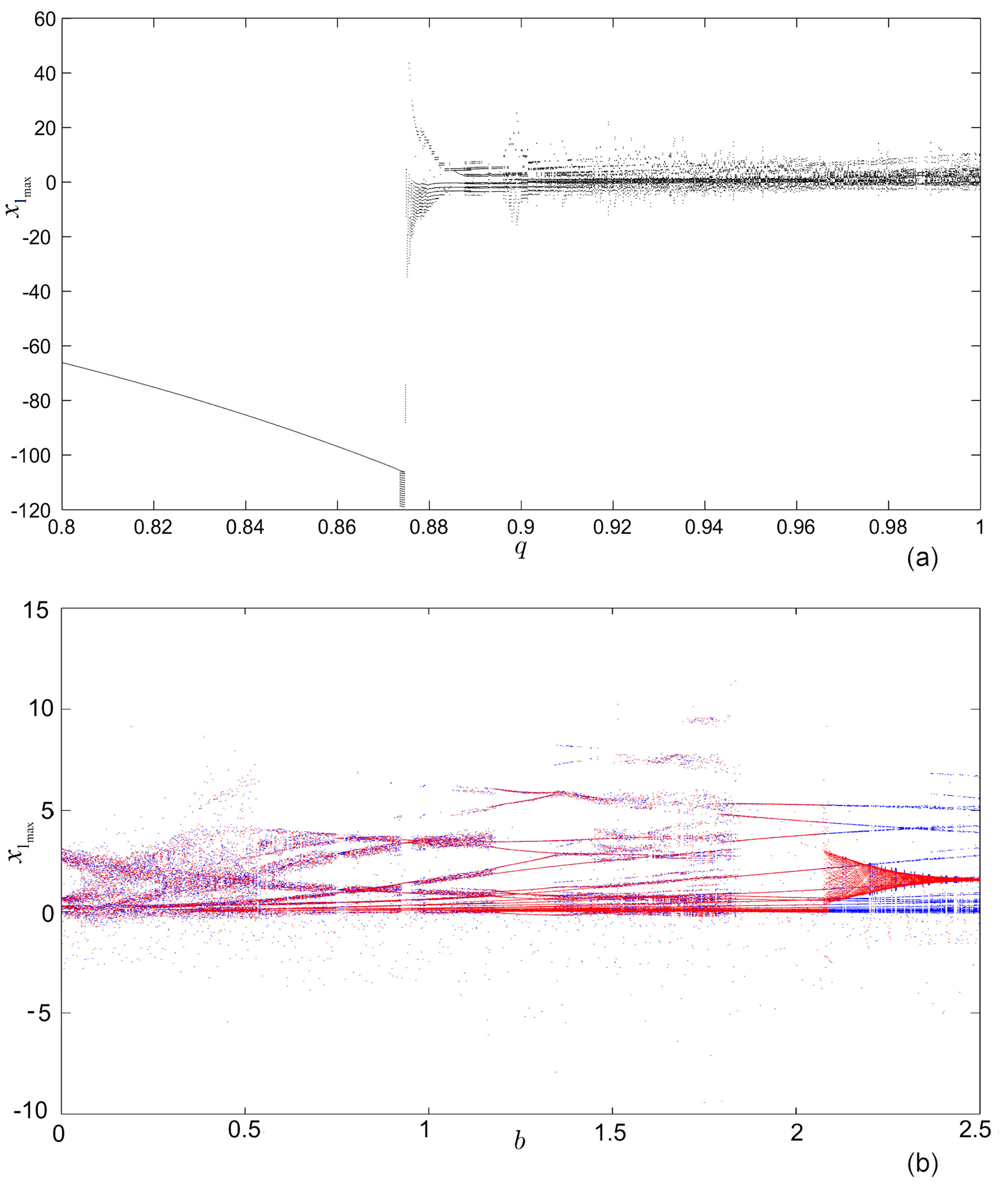}
\caption{(a) Bifurcation diagram with the fractional order $q$ as bifurcation parameter and $q=0.98$. (b) Bifurcation diagram with $b$ as bifurcation parameter. Red and blue plots indicate the multistability. }
\label{bifu}
\end{center}
\end{figure}

\begin{figure}
\begin{center}
\includegraphics[scale=0.5]{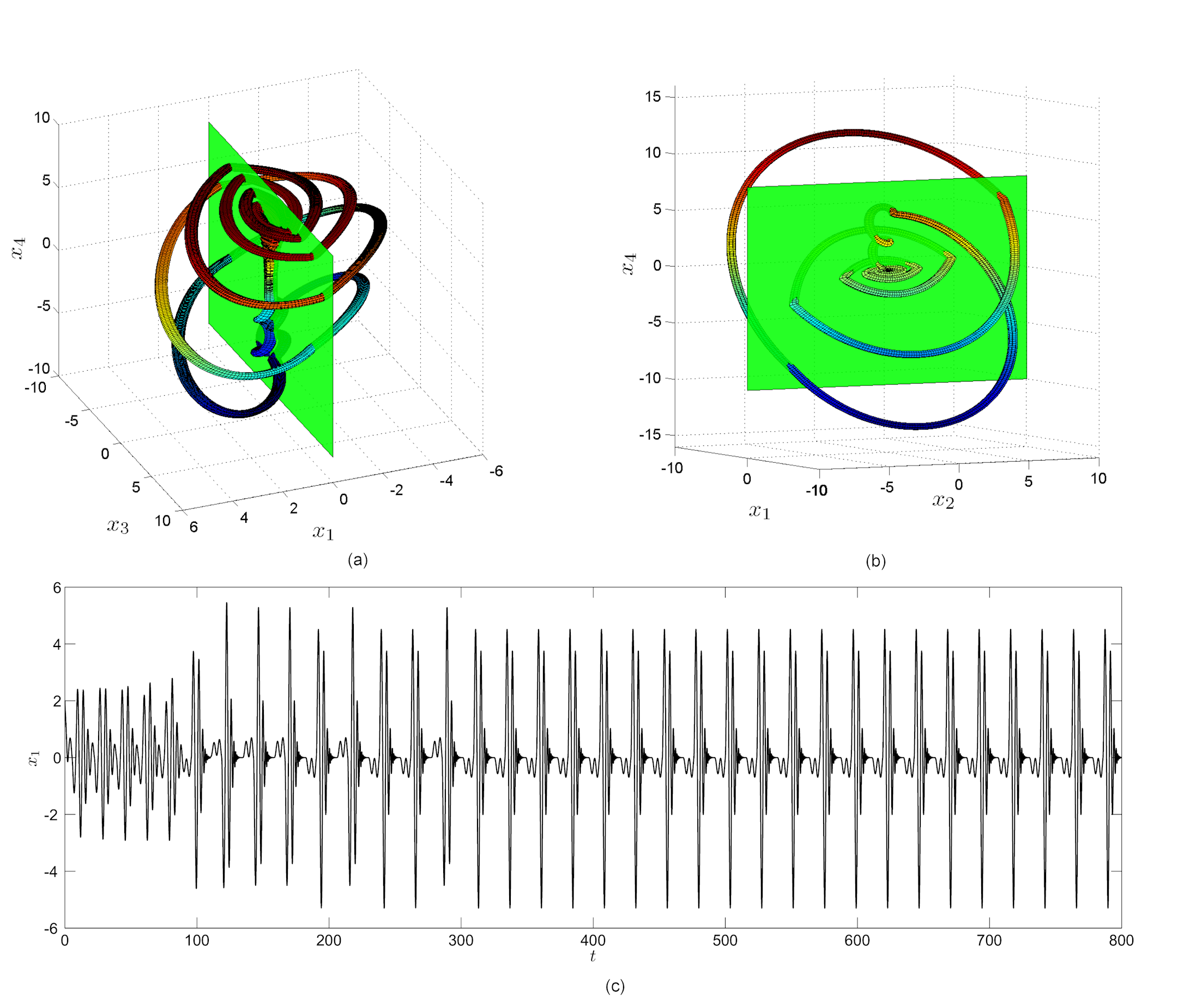}
\caption{A numeric periodic trajectory  of system \eqref{sp} for $b=2$. (a) Phase $(x_1,x_3,x_4)$ projection. (b) Phase projection $(x_1,x_2,x_4)$. The discontinuity plane $x_1=0$ reveals the corners, typical to continuous non-smooth systems, and also the sliding phenomena along the plane $x_1$.}
\label{fig2223}
\end{center}
\end{figure}

The bifurcation diagram of the approximated system, with bifurcation parameter $b$ (Fig. \ref{bifu}), shows that there are some ranges for $b$ where the system would have stable periodic cycles (see also Fig. \ref{fig2223} for $b=2$). In reality, the following result on the periodicity of solutions of FDEs should be noted.

\begin{theorem}\label{ericus}\cite{eva,moha} (see also \cite{jun})

The fractional-order system \eqref{IVP0} cannot have any exact non-constant periodic solution.
\end{theorem}

In \cite{non_eva1} and \cite{non_eva2}, it is proved that a long-time non-constant periodic solution may have a steady-state behavior, but with $-\infty$ as the lower limit in the Caputo operator. In these cases, the underlying FDEs may have asymptotically $T$-periodic solutions for which $\lim_{t\to \infty}x(t+T)-x(t)=0$ (see also \cite{peri1}), for a certain $T>0$.

For example, consider the following simple scalar linear FDE using Caputo derivative with the lower limit at $-\infty$:
\begin{equation}\label{ein11}
D^q_{-\infty}x(t)+\beta x(t)=\gamma\cos(\Omega t+\alpha),
\end{equation}
where $\alpha,\beta,\gamma,\Omega\in \mathbb{R}$ and \cite{zho}
$$
D^q_{-\infty}x(t)=\frac{1}{\Gamma(1-q)}\int_{-\infty}^t(t-s)^{-q}x'(s)ds,\quad t\in\R.
$$
A periodic solution of \eqref{ein11} is (see the proof in Appendix \ref{periodic})

\begin{equation*}%\label{ein5}
x(t)=\frac{\gamma\left(\beta\cos(\Omega t+\alpha)+\Omega^q\cos\left(\Omega t+\alpha-\frac{\pi  q}{2}\right)\right)}{\beta^2+2\beta\Omega^q\cos\frac{\pi  q}{2}+\Omega^{2 q}}.
\end{equation*}

\noindent On the other hand, the scalar linear FDE using Caputo derivative and the lower limit at $0$ ($D_*^q$),
\begin{equation*}%\label{ein6}
D^q_{*}x(t)+\beta x(t)=\gamma\cos(\Omega t+\alpha),
\end{equation*}
has not periodic solutions (Theorem \ref{ericus}).

However, continuous dynamical systems of FO, modeled by Caputo's derivative $D_*^q$, could have non-constant periodic trajectories, if the system variables are impulsed periodically (impulsive fractional-order systems) \cite{dudux}.

In consequence, these apparently periodic motions (see e.g. Fig. \ref{fig2223}), which ressemble \emph{almost periodicity} (see e.g. \cite{almost}), should be called \emph{numerically periodic oscillations}.

\begin{figure}
\begin{center}
\includegraphics[scale=0.55]{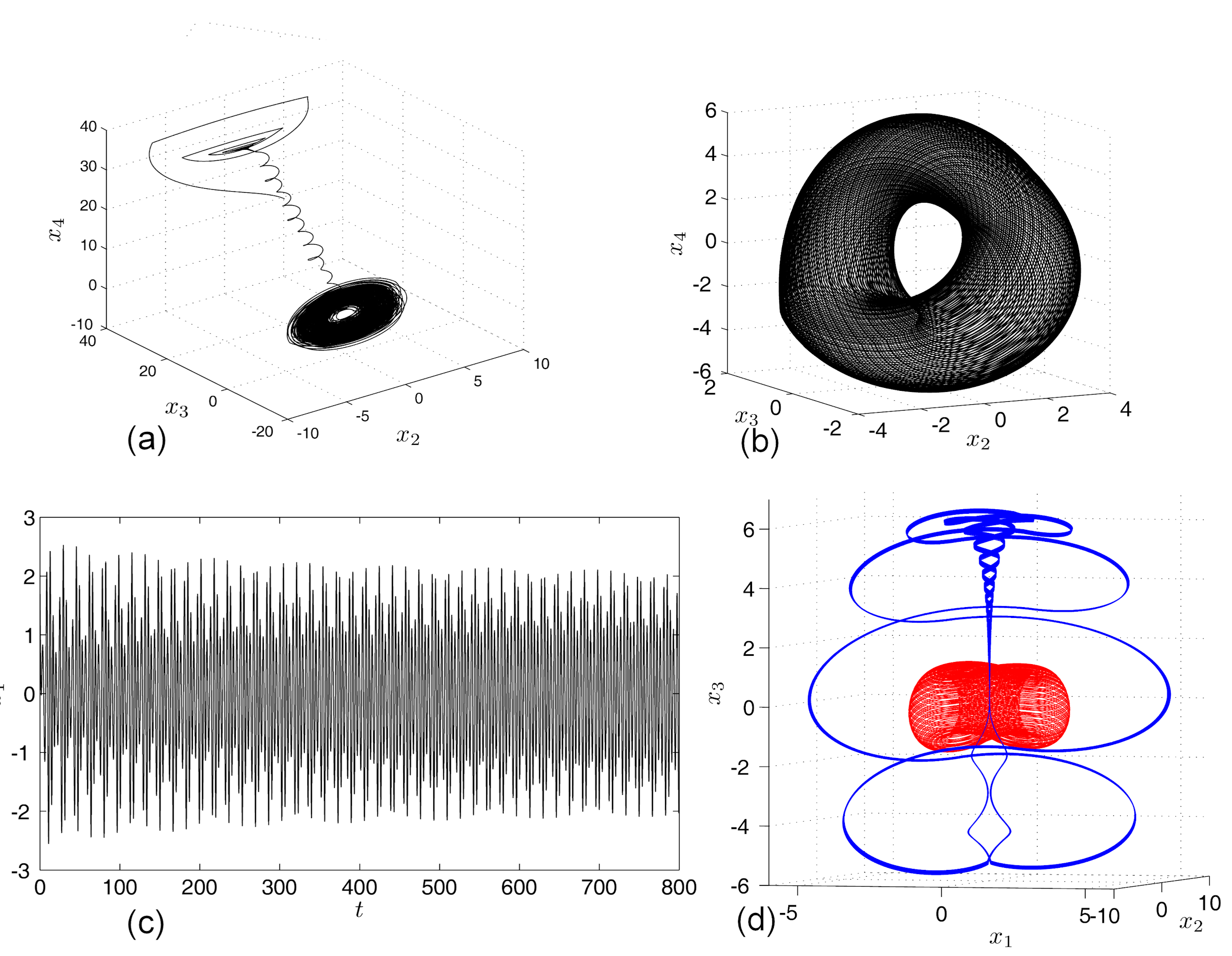}
\caption{Torus for $b=2.2$. (a) Phase space projection $(x_1,x_2,x_4)$. (b) Same torus without transient. (c) Corresponding time series. The trajectory is actually a long numerical periodic transient. (d) Two coexisting tori (red and blue respectively), for $b=2.2$, generated from different initial conditions.}
\label{torus}
\end{center}
\end{figure}

Also, the corresponding torus in the bifurcation diagram for $b>2.1$ in Fig. \ref{bifu} (b) could contain numerically  periodic oscillations (see Fig. \ref{torus}).
As remarked in the previous section, tori also present scrolls before they are reached by system trajectories (Fig. \ref{torus} (d)).

\subsection{Hidden Chaotic and hyperchaotic attractors}

From a computational point of view, it is natural
to suggest the following classification of attractors,
based on the complexity in finding basins of attraction in the phase space \cite{KuznetsovLV-2010-IFAC,LeonovKV-2011-PLA,LeonovKV-2012-PhysD,LeonovK-2013-IJBC}:
 an attractor is called a \emph{self-excited attractor}
 if its basin of attraction
 intersects with any open neighborhood of a stationary state (equilibrium);
 otherwise, it is called a \emph{hidden attractor}.

For a hidden attractor, chaotic or not, its basin of attraction is not connected with equilibria. Thus, the search and visualization of hidden attractors in the phase space could be a challenging task.

On the other hand, hidden attractors can be attractors in systems without equilibria \cite{fara} or in systems with only one stable equilibrium \cite{unul}. Hidden hyperchaotic attractors have been reported, e.g., in
\cite{BorahR-2017-HA,FengP-2017-HA,PhamVJWV-2014-HA},
and for a FO system, e.g., in \cite{VolosPZMV-2017-HA}.

In order to numerically find finite-time local Lyapunov exponents (LEs), one can locally approximate the non-smoothness, caused by the modulus function, by a quadratic polynomial $p(x)=\frac{1}{2\varepsilon}x^2+\frac{\varepsilon}{2}$ within the neighborhood $[-\varepsilon,\varepsilon]$.

\begin{figure}
\begin{center}
\includegraphics[scale=0.6]{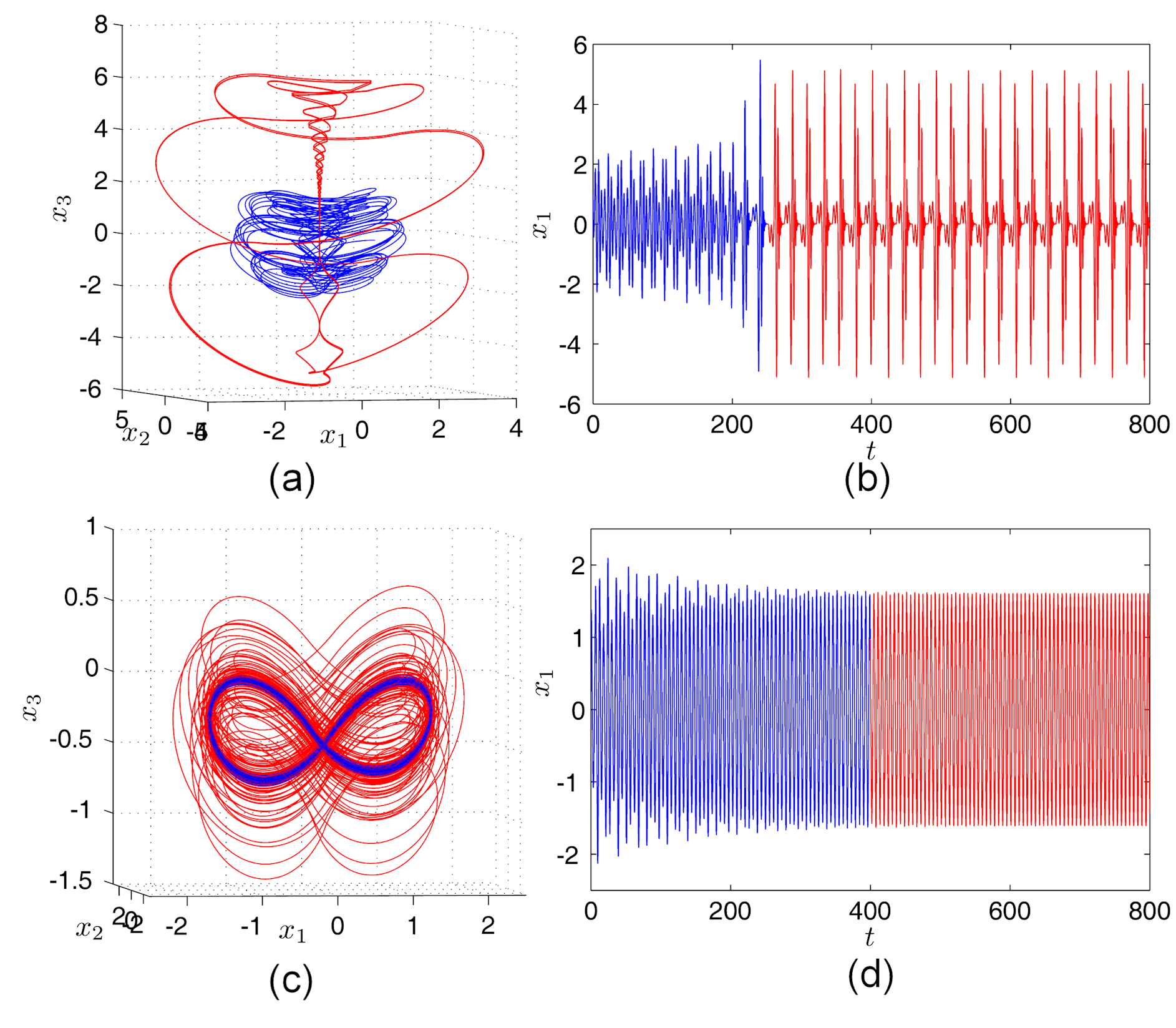}
\caption{(a) Transient hidden chaos (blue), phase projection $(x_1,x_2,x_3)$ continued by a hidden torus (red) for $q=0.9725$ and $b=1.77$. (b) Corresponding time series. (c) Transient hidden hyperchaos (blue), phase projection $(x_1,x_2,x_3)$ continued by a hidden torus (red) for $q=0.936$ and $b=1.19$. (d) Corresponding time series.}
\label{transient}
\end{center}
\end{figure}

In general, \emph{sustained chaos} is numerically indistinguishable
from \emph{transient chaos}, which can persist for a long time (see \cite{dancila1,dancila2}). For example, for $q=0.9725$ and $b=1.77$, because one of the LEs $\{0.0058,-0.0000,-0.0042,-0.0447\}$ (measured with a precision of $1E-5$ and from initial conditions $(1,2,0.0,0.1)$) is 0, the system evolves first for a relatively long time ($t\in[0,t^*]$, with $t^*\approx 220$, with hidden transient chaos, after which the trajectory is attracted by a hidden torus which, as discussed in Subsection \ref{perioada}, is characterized by a numerical periodic trajectory (see phase projections $(x_1,x_2,x_3)$ and time series in Figs. \ref{transient} (a), (b)). On the other hand, if one uses $q=0.936$ and $b=1.19$, the spectrum of the LEs (with the same initial value conditions and step size) is $\Lambda=\{0.0034,0.0022,0.0000,-0.0592)\}$, which means that the existing transient is hidden hyperchaotic, since two LEs are positive (see phase projections $(x_1,x_2,x_3)$ and the time series in Figs. \ref{transient} (c), (d)).

Because, as indicated by the bifurcation diagram (Fig. \ref{bifu}), the system presents multistability, so to obtain different hidden hyperchaotic attractors, one needs to choose the initial points in their respective basins of attraction (see, e.g. the case presented in Fig. \ref{torus} (c), where the initial conditions are $(1,2,0,1)$ and $(11,-1,0,.1)$).

For $b\in(0,{b^*})$, ${b^*}\approx1.25$, the system presents hidden hyperchaotic attractors. For example, in Fig. \ref{fig222} a hidden hyperchaotic attractor corresponding to $b=0.5$ is presented.

\begin{figure}
\begin{center}
\includegraphics[scale=0.35]{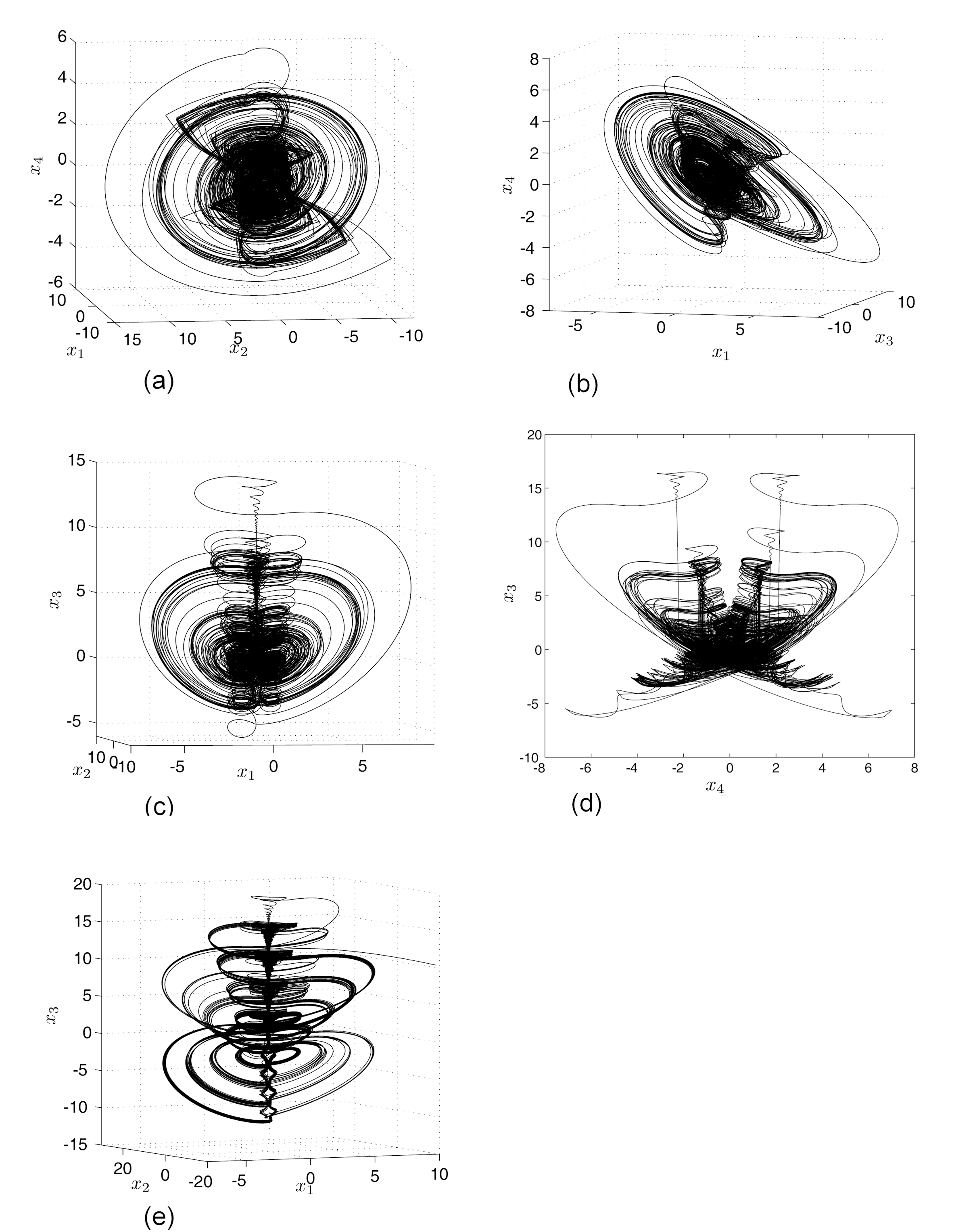}
\caption{Hidden hyperchaotic attractor corresponding to $b=0.5$. (a)-(c) 3D phase projections. (d) Plane phase projection $(x_4,x_3)$. (e) Multiple scrolls for $b=1.77$ in phase projection space $(x_1,x_2,x_3)$).}
\label{fig2}
\end{center}
\end{figure}

Also, interesting cases appear for $b\in(b_1,b_2)$, with $b_1\approx 1.45$ and $b_2\approx 1.85$. Moreover, for $b=1.75$, the trajectory presents several co-axial scrolls (Fig. \ref{fig2} (e)).

\section*{Conclusion}
In this paper, it has been shown numerically that local approximations of the discontinuities in systems \eqref{IVP0}-\eqref{f} are more useful than global approximations. Therefore, the considered system is locally continuously and smoothly approximated, for which its dynamics are numerically analyzed. Without equilibria, the system admits only hidden attractors. It has been found for what values of the fractional order $q$ and parameter $b$, the system admits a single positive finite-time Lyapunov exponent, with which the system behaves chaotically. It has also been found when the system admits two positive time-finite local Lyapunov exponents, with which the system behaves hyperchatically. Finally, it has been shown that the system cannot have exact periodic oscillations, therefore for correctness the seeming oscillations are referred to as numerically periodic oscillations.

%\numberwithin{equation}{section}
%\appendix \textbf{Appendix A}
\counterwithin{theorem}{section}
\counterwithin{definition}{section}
\counterwithin{remark}{section}
\counterwithin{equation}{section}
%\section{Appendix section}

\vspace{3mm}

\textbf{Acknowledgments}

\vspace{3mm}

M.-F. Danca, N. Kuznetsov and G. Chen are supported by the Russian Science Foundation (project 14-21-00041). M. Feckan is also supported in part by the Slovak Research and Development Agency under the Contract No. APVV-14-0378 and by the Slovak Grant Agency VEGA Nos. 2/0153/16 and 1/0078/17.

\begin{appendices}

\section{Basic notions and results}

Because the set-valued property of $F$ in (\ref{IVP1}) is generated by $S_i$, which are real functions, the notions and results presented here are considered in $\mathbb{R}$, for the case of $n=1$, but they are also valid in the general cases of $n>1$.

The graph of a set-valued function $F$ is defined as follows:
\begin{equation*}
Graph(F):=\{(x,y)\in \mathbb{R}\times \mathbb{R}, ~y\in F(x)\}.
\end{equation*}

\begin{remark}\label{rem}Due to the symmetric interpretation of a set-valued function as a graph (see e.g. \cite{aub1}), a set-valued function satisfies a property if and only if its graph satisfies it. For instance, a set-valued function is closed or convex if and only if its graph is closed or convex.
\end{remark}

\begin{definition}\label{usc}
A set-valued function $F$ is upper semicontinous (u.s.c.) at $x^0\in \mathbb{R}$ if, for any open set $B$ containing $F(x^0)$, there exists a neighborhood $A$ of $x^0$ such that $F(A)\in B$.
\end{definition}

\noindent $F$ is u.s.c. if it is so at every $x^0\in \mathbb{R}$, which means that the graph of $F$ is closed.

\begin{definition} \label{gensol} A \emph{generalized solution} to (\ref{DIgen}) is an absolutely continuous function $x:[0,T]\rightarrow\mathbb{R}$, satisfying (\ref{DIgen}) for a.a. $t\in[0,T]$.
\end{definition}

\begin{definition}\label{app}
A single-valued function $h:\mathbb{R}\rightarrow \mathbb{R}$ is called an \emph{approximation} (\emph{selection}) of the set-valued function $F$, if
\begin{equation*}
h(x)\in F(x),~~\forall	x\in \mathbb{R}.
\end{equation*}
\end{definition}

\noindent Generally, a set-valued function admits (infinitely) many approximations.

 \begin{theorem}\textbf{Cellina's Theorem} (\cite{aub1} p. 84 and \cite{aub2} p. 358)\label{theo}
 Let $F:\mathbb{R}\rightrightarrows \mathbb{R}$ have convex values $F(x)$, $x\in X$. Then, for every $\varepsilon>0$, there exists a single-valued continuous $\varepsilon$-approximation of $F$.
 \end{theorem}

\noindent See Fig. \ref{fig222} (c) for the case of the $\Sgn$ function.

\begin{theorem}\label{wei}\textbf{Weiestrass Approximation Theorem }Suppose $f$ is a continuous real-valued function defined on the real interval $[a, b]$. Then, for every $\varepsilon > 0$, there exists a polynomial $p$ such that for all $x\in [a, b]$, $|f(x) - p(x)| < \varepsilon$.
\end{theorem}

\section{Explicit solutions of IVP}\label{apB}

Because the considered system \eqref{sp} is actually PWL in each of the open half spaces $\Omega_\pm$, one can find explicit solutions (solutions existence is ensured by Theorem \ref{tui}).

Indeed, \eqref{sp} can be written as
\begin{equation}\label{e1}
D^q_*x=M_{\pm}x+m,\quad x\in\Omega_\pm,
\end{equation}
where
$$
M_{\pm}=\begin{pmatrix} -1 & 1 & 0 & 0 \\
                 0 & 0 & \mp1 & 1 \\
                 \pm1 & 0 & 0 & 0 \\
                 0 & -b & 0 & 0
\end{pmatrix},
\quad
m=-ae_3,\quad e_3=\begin{pmatrix} 0\\ 0\\ 1\\ 0
\end{pmatrix},
$$
in which \eqref{sp} is piecewise affine, so $M_+x+m=g(x)+A(x)s(x)$ for $x\in\Omega_+$ and $M_-x+m=g(x)+A(x)s(x)$ for $x\in\Omega_-$.  Then, it follows from \cite{LP} that the solution of \eqref{e1} with $x(0)=x_0$ on each $\Omega_\pm$ is given by
\begin{equation}\label{e2}
x(t)=E_q(t^qM_{\pm})x_0-a\int_0^t(t-s)^{q-1}E_{q,q}((t-s)^qM_{\pm})e_3ds,
\end{equation}
where the Mittag-Leffler matrix functions $E_{\alpha}(M_{\pm})$ and $E_{\alpha,\beta}(M_{\pm})$ are defined as \cite[p. 56]{GKMR}
$$
E_{\alpha,\beta}(M_{\pm})=\sum_{k=0}^\infty\frac{M_{\pm}^k}{\Gamma(\alpha k+\beta)},\quad E_\alpha(M_{\pm})=E_{\alpha,1}(M_{\pm}).
$$
Next, using \cite[formula (4.4.4)]{GKMR}, it follows from \eqref{e2} that
\begin{equation}\label{e3}
\begin{gathered}
x(t)=E_q(t^qM_{\pm})x_0-a\int_0^t(t-s)^{q-1}E_{q,q}((t-s)^qM_{\pm})e_3ds\\
=E_q(t^qM_{\pm})x_0-a\int_0^ts^{q-1}E_{q,q}(s^qM_{\pm})e_3ds\\
=E_q(t^qM_{\pm})x_0-at^qE_{q,q+1}(t^qM_{\pm})e_3,
\end{gathered}
\end{equation}
and
$$
\begin{gathered}
e_3=\begin{pmatrix} 0\\ 0\\ 1\\ 0
\end{pmatrix} .
\end{gathered}
$$

The last formula of \eqref{e3} gives explicit solutions of \eqref{DIfinal} on each $\Omega_\pm$, respectively.

\section{Periodic solutions for Caputo derivative with lower limit at $-\infty$}\label{periodic}
Consider the following simple scalar linear FDE with Caputo derivative with the lower limit at $-\infty$ as
\begin{equation}\label{ein1}
D^q_{-\infty}x(t)+\beta x(t)=\gamma\cos(\Omega t+\alpha),
\end{equation}
where $\alpha,\beta,\gamma,\Omega\in \mathbb{R}$. Note that \cite{zho}
$$
D^q_{-\infty}x(t)=\frac{1}{\Gamma(1-q)}\int_{-\infty}^t(t-s)^{-q}x'(s)ds,\quad t\in\R.
$$
To find a solution of \eqref{ein1} in the form of
\begin{equation}\label{ein2}
x(t)=A\cos\Omega t+B\sin\Omega t.
\end{equation}
one can use formulas \cite[(24), (26)]{non_eva2} to derive
\begin{equation}\label{ein3}
\begin{gathered}
D^q_{-\infty}x(t)=\\
\Big(A\Omega^q\cos\frac{\pi q}{2}+B\Omega^q\sin\frac{\pi q}{2}\Big)\cos\Omega t+\Big(B\Omega^q\cos\frac{\pi q}{2}-A\Omega^q\sin\frac{\pi q}{2}\Big)\sin\Omega t.
\end{gathered}
\end{equation}
Inserting \eqref{ein3} into \eqref{ein1} yields
$$
\begin{gathered}
A\Omega^q\cos\frac{\pi q}{2}+B\Omega^q\sin\frac{\pi q}{2}+\beta A=\gamma\cos\alpha,\\
B\Omega^q\cos\frac{\pi q}{2}-A\Omega^q\sin\frac{\pi q}{2}+\beta B=-\gamma\sin\alpha,
\end{gathered}
$$
which has a solution
\begin{equation}\label{ein4}
\begin{gathered}
A=\frac{\gamma\left(\beta\cos\alpha+\Omega^q\cos\left(\alpha-\frac{\pi  q}{2}\right)\right)}{\beta^2+2\beta\Omega^q\cos\frac{\pi  q}{2}+\Omega^{2q}},\\
B=-\frac{\gamma\left(\beta\sin\alpha+\Omega^q\sin\left(\alpha-\frac{\pi  q}{2}\right)\right)}{\beta^2+2\beta\Omega^q\cos\frac{\pi  q}{2}+\Omega^{2 q}}.
\end{gathered}
\end{equation}
Inserting \eqref{ein4} into \eqref{ein2}, one obtains
\begin{equation*}\label{ein5}
x(t)=\frac{\gamma\left(\beta\cos(\Omega t+\alpha)+\Omega^q\cos\left(\Omega t+\alpha-\frac{\pi  q}{2}\right)\right)}{\beta^2+2\beta\Omega^q\cos\frac{\pi  q}{2}+\Omega^{2 q}}.
\end{equation*}

\end{appendices}

%%%%% CLEAR DOUBLE PAGE!
%\newpage{\pagestyle{empty}\cleardoublepage}

%\bibliography{ref2}
%\bibliographystyle{plain}

\newpage

\newpage

%\bibliographystyle{elsarticle-num}
%\bibliography{ha}

\end{document}